\documentclass[a4paper,10pt,balance]{IEEEtran}
\usepackage{graphicx}
\usepackage{pstricks}
\usepackage{amssymb}
\usepackage{amsmath}
\usepackage{amsthm}
\usepackage{epstopdf}
\usepackage{mathptmx}
\usepackage{balance}

\DeclareMathAlphabet{\mathcal}{OMS}{cmsy}{m}{n}

\usepackage{xcolor}

\title{A Novel Time-Domain Perspective of the CPC Power Theory: Single-Phase Systems}
\author{Dimitri Jeltsema, Jacob W.~van der Woude and Marek T.~Hartman\thanks{D. Jeltsema and J.W. van der Woude are with the Department of Mathematical Physics of the Delft Institute of Applied Mathematics, Delft University of Technology, Mekelweg 4, 2628 CD Delft, The Netherlands. M.T. Hartman is with the Gdynia Maritime University, ul. Morska 81-87, 81-225, Gdynia Poland. Email: \texttt{d.jeltsema@tudelft.nl}}                                            
}

\begin{document}
\maketitle

\begin{abstract}
This paper presents a novel time-domain perspective of the Currents' Physical Components (CPC) power theory for single-phase systems operating under nonsinusoidal conditions. The proposed CPC decomposition reveals some appealing physical characteristics of the load current components in terms of measurable powers that are distributed over the branches in the load network. Instrumental in the time-domain equivalent of the reactive current is the concept of Iliovici's reactive power, which is a measure of the area of the loop formed by the Lissajous figure when plotting the current against the voltage. For sinusoidal systems, Iliovici's reactive power integral is included in the IEEE Standards on power definitions. Furthermore, the reactive current is decomposed into two new currents that represent the reactive counterparts of the active and scattered current.
\end{abstract}

\section{Introduction}

The usage of alternative sources of power has caused that the problem of energy transfer optimization has become increasingly involved with nonsinusoidal signals and nonlinear loads. The power factor (PF) is used to measure the effectiveness of the transfer of energy between an electrical source and a load. It is defined as the ratio between the power consumed by a load (real or active power), denoted as $P$, and the power delivered by a source (apparent power), denoted as $S$, i.e., 
\begin{equation}\label{eq:PF}
\text{PF} := \frac{P}{S}.
\end{equation}
The active power is defined as the average of the instantaneous power and apparent power as the product of the RMS values of the load current and source voltage. The standard approach to improve the power factor is to place a lossless compensator, such as a capacitor or an inductor, parallel to the load. Conceptually, the design of the compensator typically assumes that the source is ideal, i.e., the internal (Thevenin) impedance is negligible, producing a fixed sinusoidal voltage. 

If the load is linear and time-invariant (LTI) and the source voltage is sinusoidal, the resulting stationary  current is a shifted sinusoid, and the power factor equals the cosine of the phase-shift angle between the source voltage and current. Classically, the remaining part of the power is called reactive power, and is denoted as $Q$. The relationship between the three types of power is given by 
\begin{equation}\label{eq:Sclassical}
S^2=P^2+Q^2. 
\end{equation}
Any improvement of the PF is accomplished by the reduction of the absolute value of the reactive power, hence reducing the phase-shift between the current and the voltage.

\subsection{Budeanu's Power Model for Nonsinusoidal Signals}

For periodic nonsinusoidal voltages and currents, the problem of decomposing the apparent power into active and reactive components is much more involved. The active power in the nonsinusoidal case is given by
\begin{align}
P := \sum_{k=1}^\infty U_kI_k \cos(\phi_k),\label{eq:BudeanuP}
\end{align}
with $U_k$ and $I_k$ denoting the rms value of the $k$-th harmonic component of the voltage and the current, respectively, and $\phi_k$ denoting the angle between the $k$-th harmonic component of the voltage and the current. Inspired by the form of the active power, Budeanu \cite{Budeanu1927} defined reactive power as
\begin{align}
Q_B &:= \sum_{k=1}^\infty U_kI_k \sin(\phi_k),\label{eq:BudeanuQ}
\end{align}
and also observed that for nonsinusoidal voltages and currents the quadratic sum of the active and reactive power is not equal to the apparent power. Given this difference, a concept called distortion power $D_B^2 := S^2 -P^2 -Q_B^2$ was introducted.

\begin{figure}[t]
\begin{center}
\includegraphics[width=0.3\textwidth]{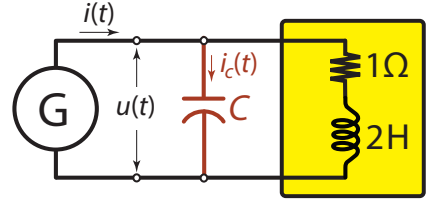}
\caption{RL circuit: uncompensated ($C=0$) and compensated ($C>0$).}
\label{fig:RLcircuit}
\end{center}
\end{figure}

\subsection{Motivation for a New Theory}\label{subsec:motivation}

Consider the RL circuit depicted in Fig. \ref{fig:RLcircuit} for $C=0$. Given a sinusoidal source voltage of the form $u(t)=10\sqrt{2}\cos{(t)}$, the associated current is given by 
\begin{equation}\label{eq:i-RLsinus}
i(t)=|Y(j)|10\sqrt{2}\cos\big(t+\phi\big),
\end{equation}
where $|Y(j)|=\frac{1}{5}\sqrt{5}$ denotes the magnitude of the load admittance ${Y}(j)=1/(1+j2)$ and $\phi=-\arctan(2)$ the associated phase shift. In this case, the active power is computed as $P=20$ W and the reactive power equals $Q=40$ VAr. Hence, the apparent power equals $S=20\sqrt{5}$ VA and the power factor equals $\text{PF}=0.447$. If a shunt capacitor is placed as a compensator, then it is clear that a capacitor of $C=0.4$ F is necessary to compensate the effect of the inductance and drive the power factor to unity. 

If the source is replaced by a nonsinusoidal voltage $u(t)=10\sqrt{2}\cos{(t)}+5\sqrt{2}\cos{(5t)}$, the associated current becomes
\begin{equation}\label{eq:i-RLnonsinus}\textstyle
i(t)=\frac{10\sqrt{2}}{\sqrt{5}}\cos\big(t+\phi_1\big)+\frac{5\sqrt{2}}{\sqrt{101}}\cos\big(t+\phi_5\big),
\end{equation}
where $\phi_1=-\arctan(2)$ and $\phi_5=-\arctan(10)$.
Using (\ref{eq:BudeanuP}) and (\ref{eq:BudeanuQ}) and $S=||u|| ||i||$, where $\|\cdot\|$ denotes the root-mean-square (rms) value (see (\ref{eq:rms})), the various power quantities and the power factor for the uncompensated circuit according to Budeanu's power model are presented in Table \ref{tab:example1}. 

As shown in the third column, although the addition of a shunt capacitor completely compensates the reactive power, the power factor is worse than in the uncompensated case. Hence the compensation of Budeanu's reactive power alone may be useless for power factor improvement as the influence of the distortion power, that only appears in nonsinusoidal periodic signals, makes this compensation far from optimal. Hence, it can be concluded that Budeanu's concept of reactive and distortion power does not exhibit any useful attributes that could be related to the physical power phenomena in the load network. Consequently, these quantities do not provide any useful information for the design of optimal compensating circuits. Although this deficiency is widely documented in the literature, see e.g., \cite{Czarnecki1987}, \cite{Shepherd1979}, 
\cite{Czarnecki2000}, \cite{EmanuelBook}, \cite{Garcia2007}, and acknowledged by the latest IEEE Standard on Power Definitions \cite{IEEE2010}, Budeanu's power model had a widespread influence for decades.

\begin{table}[t]
\caption{Power and compensation based on Budeanu's power model.}
\label{tab:example1}
\begin{center}
\begin{tabular}{|c||c|c|c|}
\hline\hline
Quantity & Uncompensated & Compensated & Unit\\
\hline
$C$ &  0       & $0.189$ & F\\
\hline
$P$ & $20.248$  & $20.248$ & W\\
\hline
$Q_B$ & $42.475$  & $0$ & VAr\\
\hline
$D_B$ & $17.800$  & $53.654$ & VA\\
\hline
$S$ & $50.309$  & $57.347$ & VA\\
\hline 
PF  & $0.403$   & $\red{0.353}$ & ---\\
\hline\hline
\end{tabular}
\end{center}
\end{table}

\subsection{Power Theory Development}

Starting from the work of Budeanu \cite{Budeanu1927}, many authors have aimed to improve the concept of reactive power in the most general case; see e.g., \cite{Czarnecki2008,EmanuelBook}, and the references therein. Most of these contributions aim at decompositions of the load current into physical meaningful orthogonal quantities. The majority of contributions use either pure frequency-domain techniques or a combination of time-domain and frequency-domain (hybrid) techniques. Time-domain techniques thus far either circumvented \cite{Garcia2007,Tano2012}, are not widely known \cite{Tenti}, or have failed to properly reveal the various power phenomena \cite{Czarnecki1987,Czarnecki2000,Czarnecki2009}. One of the most detailed work to date appears to be that of Czarnecki \cite{Czarnecki2008}, and is commonly known as the Currents' Physical Components (CPC) method. In essence, for linear time-invariant (LTI) loads, the CPC method decomposes the current into three components associated to three distinctive physical phenomena in the load: permanent energy conversion, change of load conductance with harmonic order, and phase-shift between the voltage and current harmonics.

\subsection{Contribution} 

The original CPC method uses techniques from both the time-domain and the frequency domain, and can therefore be considered as a hybrid approach. In contrast to the arguments given in, e.g, \cite{Czarnecki2008,Czarnecki2009}, in favour of the necessity of the use of frequency-domain techniques in power theory development, the results in this paper show that in the time-domain the CPC decomposition of \cite{Czarnecki2008} naturally follows by the notion of a conductance operator and a susceptance operator. These operators, in turn, provide time-domain expressions for the load conductance and the load susceptance that are the foundation of the original CPC methodology. Of key importance for the physical origin of the reactive current is the concept of Iliovici's reactive power, a concept that was introduced in \cite{Iliovici1925}---even two years before the work of \cite{Budeanu1927}. Interestingly, for sinusoidal systems, Iliovici's reactive power integral appears in the IEEE standards on power definitions \cite{IEEE2000,IEEE2010}. 

The remainder of the paper is organized as follows. Section \ref{sec:Active-Nonactive} reviews the essence of power decomposition, including the ideas of Fryze \cite{Fryze1932}, and Section \ref{sec:CPC} is included for ease of reference to the original hybrid CPC decomposition for LTI loads. The time-domain equivalent of the CPC decomposition based on the notion of the conductance and susceptance operators is presented in Section \ref{sec:timedomainpower}. Section \ref{sec:overall} proposes to further decompose the reactive current into two new currents that represent the reactive counterparts of the active and scattered current. In Section \ref{sec:RLcircuit-time}, the proposed time-domain CPC decomposition is illustrated on the RL circuit of Fig.~\ref{fig:RLcircuit}. Finally, in Section \ref{sec:conclusions}, we conclude the paper with a discussion and some directions for future research. 


\medskip

\noindent{\bfseries Notation:} Given a square integrable $T$-periodic signal $f(t)$, we define by
\begin{equation}\label{eq:rms}
||f(t)|| := \sqrt{\frac{1}{T}\int\limits_0^T f^2(t)dt}, 
\end{equation}
the root-mean-square (rms) value of $f(t)$. Voltages are represented in volts [V] and current are represented in Amp\`ere [A]. However, these units will be omitted in the text. 

\section{Active and Non-Active Power}\label{sec:Active-Nonactive}

The mathematical mechanism behind the definition of the power factor is the Cauchy-Schwarz inequality. Given any $T$-periodic square integrable source voltage $u(t)$ and the associated source current $i(t)$, this inequality takes the form\footnote{For real vectors, (\ref{eq:CS}) was proven by Cauchy (1821), and Bunyakovsky (a student of Cauchy around 1859) noted that by taking limits one can obtain an integral form of Cauchy's inequality. The general result for an inner product space was obtained by Schwarz (1885). For that reason, the integral version is sometimes referred to as the Cauchy-Bunyakovsky-Schwarz inequality.}
\begin{equation}\label{eq:CS}
||u(t)|| ||i(t)|| \geq  \frac{1}{T}\int\limits_0^T |u(t)i(t)|dt.
\end{equation}
The inequality (\ref{eq:CS}) corresponds with (\ref{eq:PF}) by identifying the left-hand side with the apparent power $S$ and the right-hand side with the absolute value of the active power $P$. The inequality (\ref{eq:CS}) holds with equality if, and only if, the source current is directly proportional (collinear) to the source voltage, and, consequently, the power factor PF equals one.  If the power factor is not equal to one, the residual of (\ref{eq:CS}) is given by (see \cite{JeltsemaECC2014} for a derivation) 
\begin{equation}\label{eq:Q^2}
\begin{aligned}
||u(t)||^2 ||i(t)||^2 &- \left( \frac{1}{T}\int\limits_0^T u(t)i(t)dt \right)^2\\
 & = \frac{1}{2T^2}\int\limits_0^T\int\limits_0^T \big(u(s)i(t) - u(t)i(s)\big)^2dsdt.
\end{aligned}
\end{equation}
The square-root of the right-hand side double integral is commonly referred to as the non-active (or useless) power and coincides with the reactive power defined by Fryze \cite{Fryze1932}. This is easily seen by decomposing the source current into two orthogonal components 
\begin{equation*}
i(t) = i_a(t) + i_F(t),
\end{equation*}
where the \emph{active} current $i_a(t)$ is the current resulting from an orthogonal projection in the direction of the source voltage
\begin{equation}\label{eq:ia}
i_a(t) = \frac{P}{||u(t)||^2}u(t),
\end{equation}
and is thus directly proportional to the source voltage, whereas the non-active (or useless) current $i_F(t)$ is orthogonal to the source voltage. As a result, we have
\begin{equation*}
||i(t)||^2 = ||i_a(t)||^2 + ||i_F(t)||^2,
\end{equation*} 
and thus that
\begin{equation*}
||u(t)||^2 ||i(t)||^2 = ||u(t)||^2 ||i_a(t)||^2 + ||u(t)||^2 ||i_F(t)||^2.
\end{equation*}

As any directly proportional components in the right-hand side of (\ref{eq:Q^2}) vanish, the non-active power reduces to
\begin{equation}\label{eq:Q^2Fryze}
Q_F^2 = \frac{1}{2T^2}\int\limits_0^T\int\limits_0^T \big(u(s)i_F(t) - u(t)i_F(s)\big)^2dsdt, 
\end{equation}
or, equivalently, $Q_F=||u(t)||||i_F(t)||$. Hence, the active, non-active, and apparent power are related via $S^2=P^2+Q_F^2$.

In sinusoidal systems, $Q_F$ is directly related to the compensating equipment ratings which can improve the power factor to unity. However, as is pointed out in \cite{Czarnecki1997}, it does not provide any (direct) useful information for compensator design in the non-sinusoidal case. For that reason, a further decomposition of the source current was proposed by decomposing $i_F(t)$ into two orthogonal quantities. This decomposition, known as the Currents' Physical Components (CPC), is briefly outlined in the next subsection.

\section{Currents' Physical Components (CPC)}\label{sec:CPC}

The development of the CPC-based power theory dates back to 1984 \cite{Czarnecki1984}, with explanations of power properties in single-phase circuits driven by a non-sinusoidal voltage source of the form
\begin{equation*}
u(t) = U_0 + \sqrt{2} \, \text{Re} \left\{\sum_{n \in N} U_n e^{jn\omega t}\right\},
\end{equation*}
where $N$ is the set of harmonics present in the signal and $\omega = 2 \pi /T$ is the fundamental frequency. The basic assumption is that the load is linear time-invariant (LTI), which can be characterized by a frequency-dependent admittance of the form 
\begin{equation}\label{eq:Y}
Y(j\omega)=G(\omega)+jB(\omega).
\end{equation}
Thus, for each harmonic $n\in N$, the associated admittance can be written as $Y_n(n\omega) = G_n(n\omega) + jB_n(n\omega)$. Consequently, the load current can be expressed as 
\begin{equation}\label{eq:CPC-totalcurrent}
i(t) = Y_0U_0 + \sqrt{2} \, \text{Re} \left\{\sum_{n \in N} Y_n U_n e^{jn\omega t}\right\}.
\end{equation}
The main idea then is to decompose the latter into three orthogonal components as
\begin{equation}\label{eq:iCPC}
i(t) = i_a(t) + i_s(t) + i_r(t).
\end{equation}  
The first component is the same as the active current defined by (\ref{eq:ia}), and can therefore be written as
\begin{equation}\label{eq:CPC:active}
i_a(t) = G_e u(t) = G_e U_0 + \sqrt{2} \, \text{Re} \left\{ \sum_{n \in N} G_e U_n e^{jn\omega t}\right\},
\end{equation}
where $G_e$ denotes the equivalent conductance 
\begin{equation}\label{eq:GeFryze}
G_e:=\frac{P}{||u(t)||^2}.
\end{equation} 
The remaining components can be found by extracting the active current from the load current, i.e.,\\
\begin{equation*}
\begin{aligned}
 i_F(t) & = i(t) - i_a(t)  \\[0.8em]
& = (Y_0-G_e) U_0 + \sqrt{2} \, \text{Re} \left\{\sum_{n \in N} (Y_n - G_e) U_n e^{jn\omega t}\right\}\\
& = (Y_0-G_e) U_0 + \sqrt{2} \, \text{Re} \left\{\sum_{n \in N} (G_n + jB_n - G_e) U_n e^{jn\omega t}\right\},
\end{aligned}
\end{equation*}
and decomposing the latter into
\begin{align}
i_s(t) &= (G_0-G_e) U_0 + \sqrt{2} \, \text{Re} \left\{\sum_{n \in N} (G_n - G_e) U_n e^{jn\omega t}\right\},\label{eq:CPC:scattered}\\
i_r(t) &= \sqrt{2} \, \text{Re} \left\{\sum_{n \in N} jB_n U_n e^{jn\omega t}\right\},\label{eq:CPC:reactive}
\end{align}
referred to as the \emph{scattered} and \emph{reactive} current, respectively. 

Although these current components are not considered as physical quantities in themselves, they are associated to three distinctive physically measurable phenomena in the load:
\begin{enumerate}
\item Permanent energy conversion---active current $i_a(t)$;
\item Change of load conductance $G_n$ with harmonic order---scattered current $i_s(t)$;
\item Phase-shift between the voltage and current harmonics---reactive current $i_r(t)$.
\end{enumerate}
For more details and a proof of orthogonality, the interested reader is referred to \cite{Czarnecki2008}.

\subsection{Active, Scattered, and Reactive Power}

Since the currents in the CPC decomposition (\ref{eq:iCPC}) are mutually orthogonal, we have
\begin{equation}\label{eq:CPC}
||i(t)||^2 = ||i_a(t)||^2 + ||i_s(t)||^2 + ||i_r(t)||^2. 
\end{equation}
This means that (\ref{eq:Q^2Fryze}) can be split into a scattered and a reactive power component as
\begin{equation*}
Q_F^2 = D_s^2 + Q_r^2,
\end{equation*}
with
\begin{align*}
D_s^2 &:= \frac{1}{2T^2}\int\limits_0^T\int\limits_0^T \big(u(s)i_s(t) - u(t)i_s(s)\big)^2dsdt,\\[0.5em]
Q_r^2 &:= \frac{1}{2T^2}\int\limits_0^T\int\limits_0^T \big(u(s)i_r(t) - u(t)i_r(s)\big)^2dsdt,
\end{align*}
or, equivalently,
\begin{align*}
D_s &= ||u(t)|| ||i_s(t)||,\\[0.5em]
Q_r &= ||u(t)|| ||i_r(t)||,
\end{align*}
yielding $S^2 = P_a^2 + D_s^2 + Q_r^2$, where we set $P_a:=P$.

\subsection{RL Circuit Revisited}\label{subsec:RL-hybrid}

Consider again the RL circuit of Fig.~\ref{fig:RLcircuit}. The associated load conductance and susceptance are given by
\begin{equation*}
G(\omega) = \frac{1}{4\omega^2 + 1}, \ B(\omega) = - \frac{2\omega}{4\omega^2+1},
\end{equation*}
respectively. If $u(t)=10\sqrt{2}\cos{(t)}+5\sqrt{2}\cos{(5t)}$, then $\omega =1$ rad/sec, $N=\{1,5\}$, and the CPC decomposition
(\ref{eq:CPC:active})--(\ref{eq:CPC:reactive}) provides the currents
\begin{equation}\label{eq:RL_CPCcurrents}
\begin{aligned}
i_a(t) &= 1.620\sqrt{2}\cos(t) + 0.810\sqrt{2}\cos(5t),\\ 
i_s(t) &= 0.380\sqrt{2}\cos(t) - 0.760\sqrt{2}\cos(5t),\\
i_r(t) &= 4.000\sqrt{2}\,\sin(t) + 0.495\sqrt{2}\,\sin(5t),
\end{aligned}
\end{equation}
whereas the active, scattered, reactive, and apparent power and the power factor for the uncompensated circuit are presented in Table \ref{tab:example2}. As shown in the third column, the addition of a shunt capacitor decreases the reactive power and slightly increases the power factor. The reason why the power factor is still less than one has two reasons: 

\begin{itemize}

\item The reactive current $i_r(t)$, and hence the associated reactive power $Q_r$, cannot be completely compensated by a capacitor (or inductor) alone. The design of a more complex compensator is necessary; 

\item Even if the reactive current would be completely compensated, the power factor will still be less than one due to the presence of the scattered current $i_s(t)$, and hence the scattered power $D_s$, which cannot be compensated with a lossless shunt compensator.\footnote{It is possible to compensate the non-active power $Q_F$ (and thus the scattered power) completely using an active filter \cite{EmanuelBook}.} 

\end{itemize}

\begin{table}[t]
\caption{Power and compensation based on the CPC power model.}
\label{tab:example2}
\begin{center}
\begin{tabular}{|c||c|c|c|}
\hline\hline
Quantity & Uncompensated & Compensated & Unit\\
\hline
$C$ &  ---       & $0.072$ & F\\
\hline
$P_a$ & $20.248$  & $20.248$ & W\\
\hline
$Q_r$ & $45.063$  & $39.467$ & VAr\\
\hline
$D_s$ & $9.505$  & $9.505$ & VA\\
\hline
$S$ & $50.309$  & $45.365$ & VA\\
\hline 
PF  & $0.403$   & $0.446$ & ---\\
\hline\hline
\end{tabular}
\end{center}
\end{table}

\section{A Time-Domain Perspective of CPC}\label{sec:timedomainpower}

Although the CPC method, as outlined in Subsection \ref{sec:CPC}, presents the voltage  and the CPC currents in the time-domain, their expressions are obtained using frequency-domain techniques. Indeed, the definition of the source voltage, as is used in the CPC decomposition \cite{Czarnecki2008}, is, due to its periodicity assumption, given in terms of a Fourier series
\begin{align}
u(t) :=& \ U_0 + \sqrt{2} \, \text{Re} \left\{\sum_{n \in N} U_n e^{jn\omega t}\right\}\notag\\
      =& \ U_0 + \sqrt{2}\sum_{n \in N} U^c_n \cos(n \omega t) + \sqrt{2}\sum_{n \in N} U^s_n \sin(n \omega t)\notag\\    =& \sum_{n \in N'} u_n(t),\label{eq:source}
\end{align}
where $N'= N \cup \{0\}$ is the set of harmonics present, including the DC term, i.e., $u_0(t)=U_0$. Although this is a pure time-domain signal, its components are often, for ease of reference, called the harmonics, but we might as well refer to these components as the voltage contents. However, the active and scattered components of the current can be deduced from the real part of the complex admittance, whereas the reactive current component follows from the imaginary part. For that reason, we refer to this approach as the \emph{hybrid} CPC method. The key question to be answered in the present section is how to establish (\ref{eq:CPC}) in the time domain and how to infer the physical origins from such decomposition. 

\subsection{Time-Domain Origins of Reactive Current}\label{subsec:reactive-time}

We notice that the definition of the reactive current (\ref{eq:CPC:reactive}) only holds for a stationary periodic situation (i.e., no transients). Inspired by the notion of the admittance operator used in \cite{Garcia2007}, we can therefore replace (\ref{eq:CPC:reactive}) by its time-domain equivalent 
\begin{equation*}
i_r(t)=\sum_{n\in N}\mathfrak{B}_n\{u_n(t)\},
\end{equation*}
where $\mathfrak{B}_n\{\cdot\}:=\frac{B_n}{n\omega}\frac{d}{dt}\{\cdot\}$ may be referred to as the susceptance operator. Note that for DC terms $B_0 \equiv 0$ and thus that $\mathfrak{B}_0\{\cdot\} \equiv 0$. Denoting the $n$-th component of $i_r(t)$ by $i_{r_n}(t)$, we have that $i_{r_n}(t)=\mathfrak{B}_n\{u_n(t)\}$. Multiplication of the latter expression on both sides with $\frac{1}{n\omega}\frac{d}{dt}u_n(t)$ and averaging over the interval $[0,T]$ yields
\begin{equation*}
\frac{1}{n\omega T}\int\limits_0^T i_{r_n}(t)\frac{du_n(t)}{dt} dt =\frac{B_n}{n^2\omega^2 T}\int\limits_0^T \left(\frac{du_n(t)}{dt}\right)^2dt,
\end{equation*}
which, in turn, suggests that $B_n$ can be expressed as 
\begin{equation}\label{eq:Bn-time}
B_n = n\omega\frac{\displaystyle\frac{1}{\omega T}\int\limits_0^T i_{r_n}(t)\frac{du_n(t)}{dt} dt}{\displaystyle\frac{1}{\omega T}\int\limits_0^T \left(\frac{du_n(t)}{dt}\right)^2dt},
\end{equation}
where $\omega T=2 \pi$.
By orthogonality of $i_{r_n}(t)$ with respect to $i_a(t)$, $i_s(t)$ and all the remaining $k$ components of $i_r(t)$ with respect to $\frac{d}{dt}u_n(t)$, for all $k \neq n$, (see e.g.,~\cite{Czarnecki2008} for a proof), there holds that
\begin{equation*}
\frac{1}{\omega T}\int\limits_0^T i_{r_n}(t)\frac{du_n(t)}{dt} dt = \frac{1}{\omega T}\int\limits_0^T i(t)\frac{du_n(t)}{dt} dt. 
\end{equation*}
Hence, the reactive current (\ref{eq:CPC:reactive}) can be written as
\begin{equation}\label{eq:reactive_time}
i_r(t) = \sum_{n \in N}\frac{\displaystyle\frac{1}{\omega T}\int\limits_0^T i(t)\frac{du_n(t)}{dt} dt}{\displaystyle\frac{1}{\omega T}\int\limits_0^T \left(\frac{du_n(t)}{dt}\right)^2dt}\frac{du_n(t)}{dt}.
\end{equation}
The latter time-domain representation of the reactive current has an appealing physical interpretation. The form of the integral in the numerator, which can be alternatively expressed as
\begin{equation}\label{eq:Iliovicipower}
\frac{1}{\omega T}\int\limits_0^T i(t)\frac{du_n(t)}{dt} dt = - \frac{1}{\omega T}\int\limits_0^T u_n(t)\frac{di(t)}{dt} dt,
\end{equation}
is equivalent to the time-parametrized version of Iliovici's reactive power integral \cite{Iliovici1925} and reflects the reactive power associated to the inductors and capacitors that is generated by the $n$-th component of the source voltage (see Subsection \ref{subsec:Iliovici} for a further discussion). Thus, for LTI loads, the reactive current is build up by the powers in the energy accumulating elements of the load network. Of course, as each of these powers is generated by the pair $\{u_n(t),i_{r_n}(t)\}$, we may also conclude, in a frequency-domain parlance, that these powers, and, consequently, the load susceptance (operator), generally change with harmonic order. Note that, starting from (\ref{eq:reactive_time}), we can arrive back at (\ref{eq:CPC:reactive}) using Fourier transform and Parseval's identity \cite{Papoulis1962}. 

\subsection{Time-Domain Origins of Active and Scattered Current}

Concerning the active current used in the CPC decomposition, the time-domain approach has no new insights to provide as this current is already defined in the time-domain by
\begin{equation}\label{eq:active}
i_a(t) = \frac{\displaystyle\frac{1}{T}\int\limits_0^T u(t)i(t)dt}{\displaystyle\frac{1}{T}\int\limits_0^T u^2(t)dt}u(t).
\end{equation}

The scattered current (\ref{eq:CPC:scattered}) can also be represented fully in the time-domain via the notion of the conductance operator $\mathfrak{G}_n=G_n-G_e$. In a similar fashion as the reactive current, we arrive at a time-domain representation of $G_n$ of the form
\begin{equation}\label{eq:Gn-time}
G_n = \frac{\displaystyle\frac{1}{T}\int\limits_0^T u_n(t)i(t) dt}{\displaystyle\frac{1}{T}\int\limits_0^T u_n^2(t) dt}.
\end{equation}
Hence, the scattered current takes the form 
\begin{equation}\label{eq:scattered_time}
i_s(t) = \sum_{n\in N'} \frac{\displaystyle\frac{1}{T}\int\limits_0^T u_n(t)i(t) dt}{\displaystyle\frac{1}{T}\int\limits_0^T u_n^2(t) dt} u_n(t) - G_e u(t).
\end{equation}
As the integral in the numerator of (\ref{eq:scattered_time}) vanishes when it is taken over the inductive and capacitive branch voltages and currents, the components of the scattered current are seen to be proportional to the powers `scattered' over all the resistive branches in the load network that are generated by each of the components of the source voltage. Again, from a frequency-domain perspective, these powers, and, consequently, the load conductance (operator), generally change with harmonic order. 

\subsection{Iliovici's Reactive Power Integral}\label{subsec:Iliovici}

As briefly highlighted in Subsection \ref{subsec:reactive-time}, a key role played in the time-domain equivalent of the CPC reactive current (\ref{eq:reactive_time}) is the concept of a reactive power integral that---to the best of our knowledge---was first proposed by Iliovici \cite{Iliovici1925}, two years before Budeanu published his power theory.\footnote{In the early fifties, Millar proposed similar integrals, called content and co-content, in order to generalize Maxwell's minimal heat theorem; see \cite{JeltsemaPhD} and the references therein. The connection of content and co-content with reactive power was also observed in \cite{Furga1994}.} Originally, this reactive power integral was given in the form
\begin{equation}\label{eq:QIliovici}
Q_I := \frac{1}{\omega T}\oint\limits_\mathcal{L} u d i,
\end{equation}
and presented for sinusoidal systems. Although it was said not to posses any physical meaning, the contour integral (\ref{eq:QIliovici}) can be considered as a measure of the area of the loop formed by a so-called Lissajous figure, i.e., the contour $\mathcal{L}$, when plotting the current against the voltage---thus it presents a quantity that can be physically measured using e.g., an oscilloscope.

\begin{figure}[t]
\begin{center}
\includegraphics[width=0.5\textwidth]{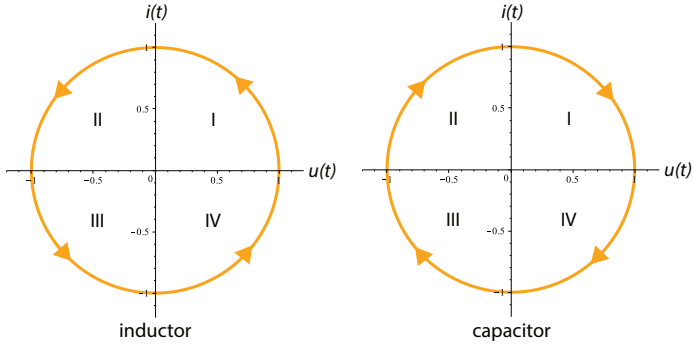}
\caption{Lissajous plot for an inductor (left) and a capacitor (right).}
\label{fig:Lissa}
\end{center}
\end{figure}

As is well-known in sinusoidal circuits, the loop area formed by plotting the current against the voltage is directly associated to the amount of phase shift produced in the load \cite{IEEE2010}. Indeed, if we consider an inductor, say with inductance $L=1$ [H], that is connected to a sinusoidal voltage source, say $u(t)=\sin(t)$, then the associated current is obtained by   
\begin{equation}\label{eq:L}
i(t) = \int\limits \frac{u(t)}{L}dt = -\cos(t),
\end{equation}
which is clearly shifted 90 degrees with respect to the voltage. It is generally known that the active power associated to the inductor is obtained by multiplying both sides of (\ref{eq:L}) by $u(t)$ and averaging over $[0,T]$, i.e., 
\begin{equation*}
P_L=\frac{1}{T}\int\limits_0^T u(t)i(t)dt=0.
\end{equation*}
Looking at the Lissajous plot of $i(t)$ versus $u(t)$ shown on the left in Fig.~\ref{fig:Lissa}, we observe that $P_L$ vanishes due to the fact that $P_L$ represents the sum of the positive powers associated to the first and the third quadrant and the negative powers associated to the second and fourth quadrant. However, a measure of the total area of the Lissajous plot is given by (\ref{eq:QIliovici}). For the inductor this measure will be non-zero since the orientation is consistently anti-clockwise. Indeed, as (\ref{eq:QIliovici}) can also be represented by its time-parameterized equivalents
\begin{equation*}
Q_I  = \frac{1}{\omega T}\int\limits_0^T u(t)\frac{di(t)}{dt}d t = -\frac{1}{\omega T}\int\limits_0^T i(t)\frac{du(t)}{dt}d t,
\end{equation*}
we compute 
\begin{equation*}
Q_I = \frac{1}{\omega T}\int\limits_0^T \sin(t)\frac{d-\cos(t)}{dt}d t =\frac{1}{2}.
\end{equation*}
Note that the sign of $Q_I$ can be used to characterize the nature of the load. It is well-known that in an inductor the current lags the voltage, which results in $Q_I > 0$ due to 
the anti-clockwise orientation of the Lissajous figure. Note that all of the above dually applies to a capacitor starting from 
\begin{equation}\label{eq:C}
i(t) = C \frac{du(t)}{dt} = \cos(t). 
\end{equation}
Indeed, applying the same voltage to a capacitor $C=1$ [F], the current equals $i(t) = \cos(t)$, and hence
\begin{equation*}
Q_I = \frac{1}{\omega T}\int\limits_0^T \sin(t)\frac{d\cos(t)}{dt}d t =-\frac{1}{2}.
\end{equation*}
Now, $Q_I < 0$, which implies that the voltage lags the current due to a clock-wise orientation of the Lissajous plot shown on the right of Fig.~\ref{fig:Lissa}. The same conclusions apply to nonlinear time-invariant inductors and capacitors. On the other hand, if we replace the inductor or the capacitor by a linear time-invariant resistor, then $u(t)=Ri(t)$ or $i(t)=Gu(t)$, and it is easy to shown that $Q_I=0$ as there is no closed contour present in the associated Lissajous plot. The same conclusion applies to all nonlinear time-invariant resistors \cite{Wyatt}. However, linear or nonlinear time-variant resistors may generate orientationally-consistent closed contours in the current-voltage plane, hence cause a phase shift between voltage and current, and, as a result, give rise to reactive power. Thus, the common belief that reactive power is due to energy accumulation only holds for (linear or nonlinear) time-invariant (TI) loads. An example of a load exhibiting reactive power without energy storage is given in \cite{Czarnecki2007}. 

\subsection{Comparison to Budeanu's Concept of Reactive Power}

Under sinusoidal conditions, Iliovici's reactive power integral (\ref{eq:QIliovici}) measures the area of the loop formed by the Lissajous plot of the current versus the voltage. The fundamental frequency $\omega$ determines the number of times the loop is cycled. To normalize the area to one cycle, the total area is divided by the fundamental frequency $\omega$. Under nonsinusoidal conditions, the reactive power associated with $n$-th harmonic should be measured relative to the fundamental harmonic as is realized by (\ref{eq:Iliovicipower}), i.e.,
\begin{equation*}
Q_{I_n}  := \frac{1}{\omega T}\int\limits_0^T u_n(t)\frac{di(t)}{dt}d t.
\end{equation*}
This means that the $n$-th harmonic produces $n$ times more loop cycles than the fundamental harmonic. Thus, the magnitudes of the reactive powers increase---relative to the reactive power associated to the fundamental harmonic---with harmonic order. Budeanu's reactive power (\ref{eq:BudeanuQ}), however, normalizes all the individual reactive powers to the fundamental frequency, i.e.,
\begin{equation*}
Q_B := \sum_{n \, \in \, N} ||u_n(t)||||i_n(t)|| \sin(\phi_n) = 
\sum_{n \, \in \, N} \frac{Q_{I_n}}{n}.
\end{equation*}
Hence, reactive powers associated to voltage harmonics with the same rms value, but different frequencies, may cancel out, and thus provide an erroneous picture of the load characteristics \cite{Czarnecki1987,IEEE2010}.

\subsection{Equivalent Susceptance}

Interestingly, Iliovici's reactive power integral (\ref{eq:QIliovici}) can be considered as the reactive equivalent of the active power. This is apparent for sinusoidal systems as (\ref{eq:QIliovici}) is precisely equal to the classical definition of reactive power \cite{IEEE2010}. For nonsinusoidal systems, however, the integral in the form (\ref{eq:QIliovici}) presents an average taken over all the voltage components (harmonics). This suggests to extend the original CPC decomposition with the notion of an equivalent load susceptance 
\begin{equation}\label{eq:Be}
B_e :
= \omega\dfrac{\displaystyle\frac{1}{\omega T}\int\limits_0^T i(t)\frac{du(t)}{dt}d t}{\displaystyle\frac{1}{\omega T}\int\limits_0^T \left(\frac{du(t)}{dt}\right)^2 d t}.
\end{equation}
Note that the concept of an equivalent susceptance is in some sense dual to the concept of equivalent conductance (\ref{eq:GeFryze}) as originally defined by Fryze \cite{Fryze1932}. The equivalent conductance represents a measure of the average conductance characteristic that is responsible for the amount of active power dissipated ($P>0)$ or delivered ($P<0$) by the load. Hence, if $G_e >0$, the load is dominantly passive, and if $G_e<0$, the load is dominantly active, respectively. The equivalent susceptance, on the other hand, represents a measure of the average susceptance and can be used to characterise the load behavior, i.e.,

\medskip

\begin{itemize}
\item $B_e=-\dfrac{\omega^2 Q_I}{\|\dot{u}(t)\|^2} < 0$ (dominantly inductive);\\[0.5em]  
\item $B_e=-\dfrac{\omega^2 Q_I}{\|\dot{u}(t)\|^2} > 0$ (dominantly capacitive), 
\end{itemize}

\medskip

\noindent where $\dot u(t) = \frac{du(t)}{dt}$.

Compensation based on placing a shunt capacitor (or inductor) that eliminates $B_e$ by setting $C=|B_e|/\omega$ (or $L=\omega/|B_e|$) reduces the average reactive power (\ref{eq:QIliovici}) to zero. Interestingly, this approach is equivalent with the compensation technique  known as energy-equalization \cite{Garcia2007}. However, as illustrated in Section \ref{sec:RLcircuit-time}, in nonsinusoidal situations this approach generally does not fully reduce the reactive power $Q_r$. 

\subsection{Equivalent Load Admittance}

The introduction of the equivalent susceptance $B_e$, together with the equivalent conductance $G_e$, naturally suggests the notion of an equivalent load admittance 
\begin{equation*}
Y_e = G_e + j B_e. 
\end{equation*}
The equivalent load admittance represents the average load behaviour and constitutes four different types of load characteristics: passive-capacitive, passive-inductive, active-capacitive, or active-inductive.   

\section{Proposed CPC Decomposition}\label{sec:overall}

\begin{figure}[t]
\begin{center}
\includegraphics[width=0.45\textwidth]{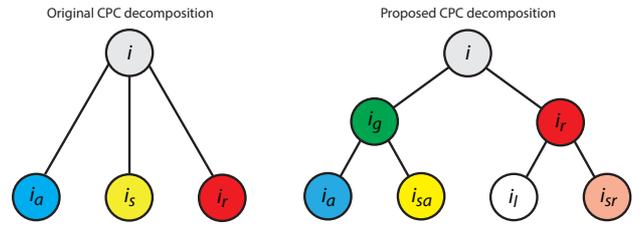}
\caption{CPC decomposition: original (left) versus proposed (right).}
\label{fig:CPCoverview}
\end{center}
\end{figure}

The CPC time-domain equivalents of the active, scattered, and reactive currents given by (\ref{eq:active}), (\ref{eq:scattered_time}), and (\ref{eq:reactive_time}), can be compactly written as
\begin{align}
i_a(t) &= \frac{1}{T} \sum_{n \in N'} \frac{u(t)}{||u(t)||^2} \int\limits_0^T u_n(t) i(t)dt,\label{eq:ia-time}\\ 
i_s(t) &= \frac{1}{T} \sum_{n \in N'} \left(\frac{u_n(t)}{||u_n(t)||^2} - \frac{u(t)}{||u(t)||^2}\right) \int\limits_0^T u_n(t) i(t)dt,\label{eq:is-time}\\
i_r(t) &= \frac{1}{T} \sum_{n \in N} \frac{\dot u_n(t)}{||\dot u_n(t)||^2}  \int\limits_0^T \dot u_n(t) i(t)dt,\label{eq:ir-time}
\end{align}
respectively, where $\dot u(t) = \frac{du(t)}{dt}$, and noting that
\begin{equation*}
\int\limits_0^T u(t) i(t)dt = \sum_{n \in N'} \int\limits_0^T u_n(t) i(t)dt.
\end{equation*}
As the latter integrals vanish when taken over the inductive and capacitive branch voltages and currents, the active and scattered power are only related to the power circulating along the resistive (conductive) branches. The total current associated to the resistive (conductive) power is determined for each voltage component by $G_n$ and is thus given by the sum of the active and scattered current \cite{HartmanCompel}, i.e.,
\begin{equation}\label{eq:ig-time}
i_g(t) = i_a(t)+i_s(t) = \frac{1}{T} \sum_{n \in N'} \frac{u_n(t)}{||u_n(t)||^2}\int\limits_0^T u_n(t) i(t)dt.
\end{equation}
This current is already implicitly present in the original CPC decomposition and it can be considered as the counterpart of the reactive power. A schematic overview of the original CPC decomposition is given in Fig.~\ref{fig:CPCoverview} (left). 

Based on the notion of the equivalent susceptance (\ref{eq:Be}), the reactive current (\ref{eq:ir-time}) can be further decomposed into a current 
\begin{equation}\label{eq:iL-time}
i_I(t) = \frac{1}{T} \sum_{n \in N} \frac{\dot u(t)}{||\dot u(t)||^2}  \int\limits_0^T \dot u_n(t) i(t)dt,
\end{equation}
which we propose to refer to as the \emph{Iliovici current}. In a similar fashion as the definition of the scattered current, the extraction of the Iliovici current (\ref{eq:iL-time}) from the reactive current (\ref{eq:ir-time}) yields
\begin{equation}\label{eq:Xcurrent}
i_{sr}(t) = \frac{1}{T} \sum_{n \in N}\left(\frac{\dot u_n(t)}{||\dot u_n(t)||^2} - \frac{\dot u(t)}{||\dot u(t)||^2} \right)\int\limits_0^T \dot u_n(t) i(t)dt.
\end{equation}
Since the integral in the numerator of (\ref{eq:Xcurrent}) vanishes when taken over the resistive branch voltages and currents, the components of (\ref{eq:Xcurrent}) are seen to be proportional to the powers `scattered' over all the inductive and capacitive branches in the load network that are generated by each of the components of the source voltage. Consequently, we propose to refer to (\ref{eq:Xcurrent}) as the \emph{scattered reactive} current. By orthogonality, the holds that
\begin{equation*}
||i_r(t)||^2=||i_I(t)||^2+||i_{sr}(t)||^2.
\end{equation*}
The associated normed Iliovici power is defined as
\begin{equation}\label{eq:Qi}
Q_i := ||u(t)|| ||i_I(t)||.
\end{equation}
Note that we intendently used a lower cast subscript `$i$' instead of `$I$' as in (\ref{eq:QIliovici}), because (\ref{eq:QIliovici}) represents the averaged reactive power, which is, in general, different from (\ref{eq:Qi}). The associated scattered reactive power equals
\begin{equation}\label{eq:Qx}
Q_s := ||u(t)|| ||i_{sr}(t)||.
\end{equation}
Overall, we have $S^2=P_a^2+D_s^2+Q_i^2+Q_s^2$ and $Q_r^2=Q_i^2+Q_s^2$. 

In order to fully `symmetrize' the CPC decomposition, we propose to refer to the original scattered current (\ref{eq:is-time}), as the \emph{scattered active} current $i_{sa}(t)$, i.e., $i_{sa}(t) \equiv i_s(t)$. The complete overview of the proposed CPC decomposition is depicted in Fig.~\ref{fig:CPCoverview} (right). From a mathematical perspective, the currents (\ref{eq:ia-time})--(\ref{eq:Xcurrent}) are nothing else than orthogonal projections (Gram-Schmidt orthogonalization in a function space) with respect to the source voltage components and their time-derivatives.

\section{RL Circuit Revisited (Closure)}\label{sec:RLcircuit-time}

\begin{figure*}[t]
\begin{center}
\includegraphics[width=0.75\textwidth]{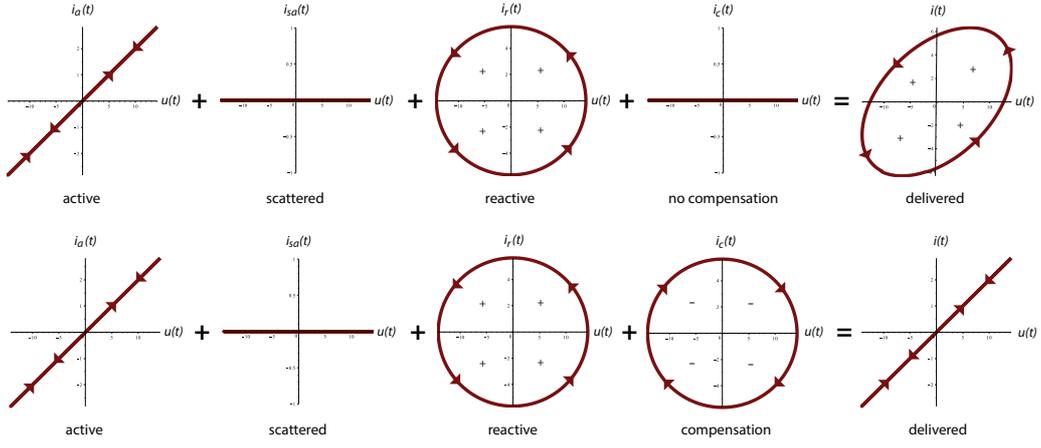}
\caption{Lissajous figures for the series RL circuit (sinusoidal): uncompensated (top) and compensated (bottom).}
\label{fig:RL_lissa_sinus}
\end{center}
\end{figure*}

\subsubsection{Sinusoidal Case}

Let us return to the series RL circuit of Fig.~\ref{fig:RLcircuit}. For a sinusoidal source voltage, it is easily checked that the computation of (\ref{eq:ia-time}) and (\ref{eq:ir-time}) provide the same active and reactive current as obtain in Subsection \ref{subsec:motivation}, respectively. Of course, in the sinusodal case, the scattered active current (\ref{eq:is-time}) equals zero. As the equivalent conductance $G_e = 0.4$, the load is dominantly passive. This is not surprising as we are dealing with a non-negative load resistor and there are no other sources present than the supply voltage. Furthermore, (\ref{eq:Be}) provides an equivalent susceptance $B_e=-0.4$, which suggests that the load is dominantly inductive. This is also not surprising since the load consists only of a resistor-inductor combination, and for that reason a natural choice to improve the power factor is to place a shunt capacitor with capacitance $C=|B_e|/\omega=0.4$ [F], which coincides with the value obtained in Subsection \ref{subsec:RL-hybrid}. 

The reactive power $Q$ equals the reactive power generated by the load inductor, i.e., $Q=Q_L$, which, in turn, equals the average reactive power $Q_L=Q_I$. In this case, it is sufficient to compensate the average reactive power since the average reactive power generated by the shunt capacitor equals $Q_C=-Q_I$. The reason why, in the sinusoidal case, it is sufficient to compensate the average reactive power is due to the fact that the Lissajous plots associated to the reactive power generated by the load inductor bears the same shape as the Lissajous plot associated to the shunt capacitor; see the Lissajous plots depicted in Fig.~\ref{fig:RL_lissa_sinus}.

\subsubsection{Nonsinusoidal Case}

\begin{figure}[t]
\begin{center}
\includegraphics[width=0.3\textwidth]{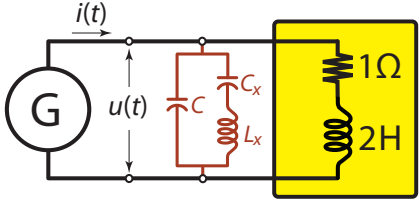}
\caption{RL circuit: Optimal compensation of reactive power.}
\label{fig:RLcomp_full}
\end{center}
\end{figure}

\begin{figure*}[t]
\begin{center}
\includegraphics[width=0.75\textwidth]{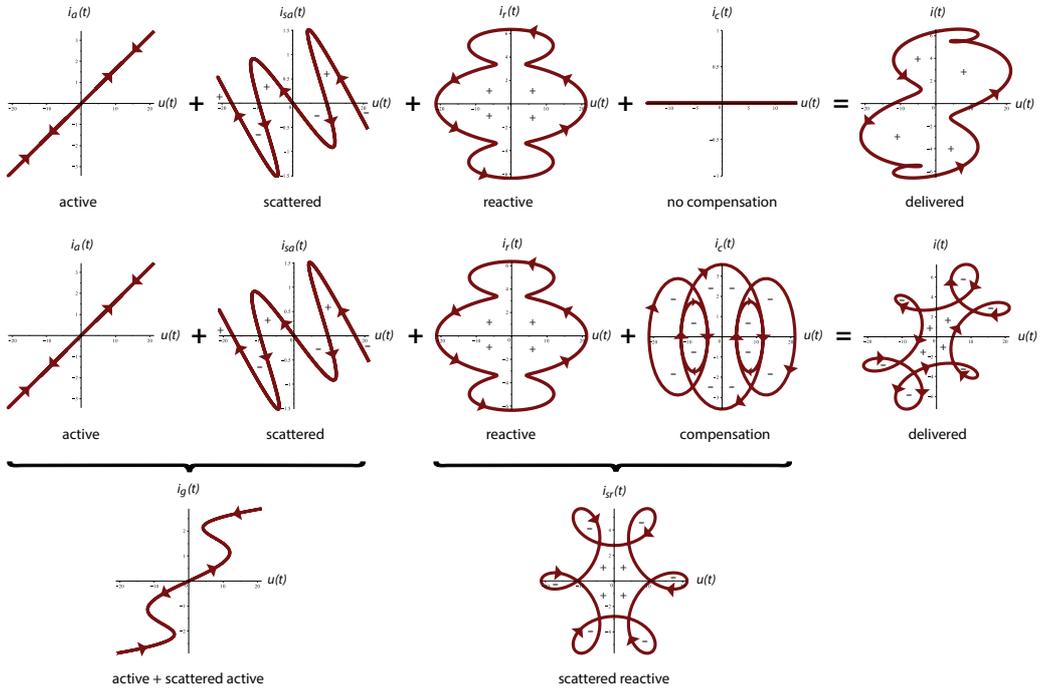}
\caption{Lissajous figures for the series RL circuit (nonsinusoidal): uncompensated (top) and compensated (bottom).}
\label{fig:RL_lissa_nonsinus}
\end{center}
\end{figure*}

Concerning the nonsinusoidal situation discussed in Subsection \ref{subsec:RL-hybrid}, it is easily checked that the computation of (\ref{eq:ia-time}), (\ref{eq:is-time}), and (\ref{eq:ir-time}) provide the same active, sctattered, and reactive current, respectively, as obtain in (\ref{eq:RL_CPCcurrents}). Now, $G_e = 0.162$ and $B_e=-0.072$, which again suggest that the load is dominantly passive and inductive. For that reason, the power factor is improved by placing a shunt capacitor with capacitance $C=|B_e|/\omega=0.072$ [F]---the same value as obtained in Subsection \ref{subsec:RL-hybrid} via minimization of the reactive power $Q_r$. The shunt capacitor eliminates the Iliovici current (\ref{eq:iL-time}) 
\begin{equation*}
i_I(t) = 0.722\sqrt{2}\sin(t)+1.806\sqrt{2}\sin(5t),
\end{equation*}
and hence the Iliovici reactive power is compensated by $Q_C=-Q_I=-52.376$ [VAr]. The associated normed power (\ref{eq:Qi}) is $Q_i=21.748$ [VAr]. As is already observed in Subsection \ref{subsec:RL-hybrid}, the reactive power seen from the generator is reduced but still not vanished. This is due to the fact that, although the shunt capacitor fully compensates the average reactive power generated by the load inductor, the shapes of the corresponding Lissajous figures depicted in Fig.~\ref{fig:RL_lissa_nonsinus} do not match. Although the average reactive power seen by the generator is zero as the Iliovici current $i_I(t)=0$), there is still a considerable amount of phase-shift as is observed from the bottom most right Lissajous plot. The remaining current is represented by the scattered reactive current (\ref{eq:Xcurrent}); see Fig.~\ref{fig:RL_lissa_nonsinus} (bottom). Indeed, 
\begin{equation*}
i_{sr}(t)=3.278\sqrt{2}\sin(t)-1.311\sqrt{2}\sin(5t),
\end{equation*}
and thus the remaining scattered reactive power equals $Q_s=39.467$ [VAr] (compare with $Q_r$ in Table \ref{tab:example2}). 

The reactive power $Q_r$ can be completely compensated by the addition of a series LC network as shown in Fig.~\ref{fig:RLcomp_full}. The values of $L_x$ and $C_x$ are computed from the set of differential equations
\begin{align*}
-L_x \frac{di_{sr}(t)}{dt} &= u(t)-u_{C_x}(t),\\
C_x \frac{du_{C_x}(t)}{dt} &= -i_{sr}(t),
\end{align*}
where $u_{C_x}(t)$ denotes the voltage across the capacitor $C_x$, and result in $L_x= 0.922$ [H] and $C_x = 0.252$ [F]. The overall results are summarized in Table \ref{tab:example3}. Note that the Lissajous plot generated by the series LC network has the same shape as the `scattered reactive' plot at the bottom-right of Fig.~\ref{fig:RL_lissa_nonsinus}, but with opposite (clockwise) orientation. The resulting load seen from the generator is characterised by the `active + scattered active' Lissajous plot on the bottom-left of Fig.~\ref{fig:RL_lissa_nonsinus}. 

Of course, in a practical situation, where the load is generally changing with time, the LC filter is replaced by an active filter. Such active filter can also be used to eliminate the scattered active current as to fully optimize the power factor.   

\begin{table}[t]
\caption{Reactive power compensation based on the proposed CPC power model.}
\label{tab:example3}
\begin{center}
\begin{tabular}{|c||c|c|c|c|}
\hline\hline
Quantity & Uncomp. & Iliovici comp. & Full comp. & Unit\\
\hline
 $C$ &  ---     & $0.072$  & $0.072$ & F\\
\hline 
$L_x$ & ---     & ---  & $0.922$ & H\\
$C_x$ & --- 	& ---  & $0.252$ & F\\
\hline
$P_a$ & $20.248$  & $20.248$ & $20.248$ & W\\
\hline
$D_s$ & $9.505$  & $9.505$ & $9.505$ & VA\\
\hline
$Q_i$ & $21.748$ & $0$ & $0$ & VAr\\
$Q_s$ & $39.467$ & $39.467$ & $0$ & VAr\\
\hline
$Q_r$ & $45.063$  & $39.467$ & $0$ & VAr\\
\hline
$S$ & $50.309$  & $45.365$ & $22.368$ & VA\\
\hline 
PF  & $0.403$   & $0.446$ & $0.905$ & ---\\
\hline\hline
\end{tabular}
\end{center}
\end{table}

\section{Concluding Remarks and Outlook}\label{sec:conclusions}

This paper presents a novel perspective of the CPC decomposition \cite{Czarnecki2008} in the time-domain, and proposes to decompose the reactive current into two additional components: an average reactive current and a scattered reactive current. Furthermore, for the sake of symmetry, the original CPC scattered current is proposed to be referred to as active scattered current. The time-domain CPC decomposition naturally follows from the notion of a conductance and susceptance operator. Of key importance is the concept of Iliovici's reactive power integral, which is a measure of the area of the loop formed by the Lissajous plot of the current against the voltage. For sinusoidal systems, this loop area is directly associated to the amount of phase-shift produced by the load \cite{IEEE2010}. For sinusoidal systems, Iliovici's integral is equivalent with the reactive power. However, for nonsinusoidal systems Iliovici's integral represents the average of the reactive powers associated to each harmonic present in the source voltage. 

\begin{figure*}[t]
\begin{center}
\includegraphics[width=0.7\textwidth]{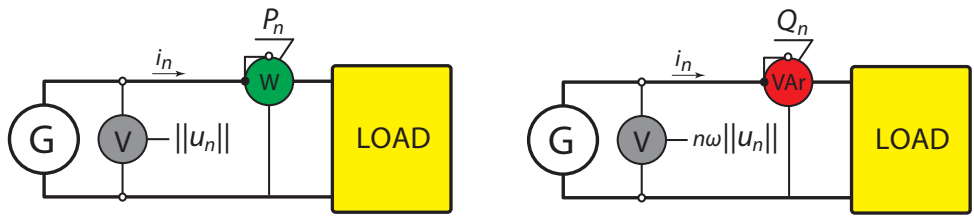}
\caption{Test circuits for measurement of the active power (left) and the reactive power (right) for the $n$-th order harmonic. Note that $||\dot{u}_n(t)||=n\omega ||u_n(t)||$.}
\label{fig:power_measure}
\end{center}
\end{figure*}

In contrast to our previous work \cite{JeltsemaECC2014}, which considers a time-domain CPC decomposition that is based on an even--odd decomposition of the load's impulse response, the results presented herein show that the origin of $i_g(t)$ (active + scattered active current) and $i_r(t)$ (Iliovici + scattered reactive current) can be related to distinctive powers distributed in the load that can be observed from measurements:

\begin{itemize}
\item The current $i_g(t)$ represents the equivalent portion of the currents that are proportional to the ratio of the powers that are circulating among the resistive branches of the load and the associated rms values of the voltage. These quantities can be measured as illustrated in Fig.~\ref{fig:power_measure} (left). Note that the measured $P_n$ provides the value of the integrals appearing in (\ref{eq:ig-time}), i.e.,
\begin{equation*}
\frac{1}{T}\int\limits_0^T u_n(t)i_n(t)dt = P_n.
\end{equation*}
\item The current $i_r(t)$ represents the equivalent portion of the currents that are proportional to the powers that are circulating among the inductive and capacitive branches of the load and the associated rms values of the time-derivative of the voltage. These quantities can be measured as illustrated in Fig.~\ref{fig:power_measure} (right). Note that the measured $Q_n$ provides the value of the Iliovici reactive power integrals appearing in (\ref{eq:ir-time}), i.e.,
\begin{equation*}
-\frac{1}{2 \pi}\int\limits_0^T \dot{u}_n(t)i_n(t)dt = Q_n.
\end{equation*}
 
\end{itemize} 

Based on the earlier work of the first author in \cite{JeltsemaPhD}, it can be shown that the power integrals appearing in the numerators of (\ref{eq:Bn-time}) and (\ref{eq:Gn-time}) can be obtained by summing the powers over the individual elements (resistors, inductors, and capacitors) in the load. Thus, these powers satisfy Tellegen's theorem. This naturally leads to the suggestion that the CPC decomposition may be converted into a conservative power theory, which is not surprising as the exposition of Sections \ref{sec:timedomainpower}--\ref{sec:overall} appears to be closely related to the Conservative Power Theory (CPT) developed by Tenti et al. \cite{Tenti}. 

Concerning power factor optimization using lossless shunt compensators, the main conclusion drawn from the proposed CPC decomposition can be summarized as follows:

\medskip

\begin{center}
\fbox{
\begin{minipage}[c]{0.45\textwidth}
The reactive power of the load can be fully compensated by adding a shunt compensator that generates the same Lissajous plot, but with opposite orientation.
\end{minipage}}
\end{center}

\medskip

\noindent For nonsinusoidal systems, this means that compensation using a single lossless shunt element, like a (bank of) capacitor(s), or an inductor, only compensates the average reactive power. More complicated shunt architectures are necessary to fully compensate the (scattered) reactive power.  

The extension to nonlinear and time-varying loads is currently studied. The results, as well as the extension of the time-domain CPC decomposition to poly-phase systems, will be reported elsewhere.

\end{document}